\documentclass{article}

\usepackage{graphicx}      

\usepackage{graphicx}

\usepackage{savesym}
\savesymbol{AND}
\usepackage{algorithm, algorithmic, setspace}
\usepackage{graphicx} 
\usepackage{amsmath} 
\usepackage{amssymb}  
\usepackage{amsfonts}
\usepackage{color}
\usepackage[english]{babel}

\newcommand{\argmin}[1]{\underset{#1}{\mathrm{argmin\,}}}

\newcommand{\xtrue}{\widetilde{x}}

\newcommand*{\QED}{\hfill\ensuremath{\square}}%

\newcommand{\xstar}{x^{\star}}
\newcommand{\R}{\mathbb{R}}

\newcommand{\fun}{\mathcal{F}}
\newcommand{\sun}{\mathcal{S}}
\newcommand{\soft}{\mathrm{S}}
\newcommand{\gun}{\mathcal{G}}
\newcommand{\run}{\mathcal{R}_{\lambda}}
\newcommand{\runinitial}{\mathcal{R}_{\alpha}}

\newcommand{\supp}{\mathcal{S}}
\newcommand{\prox}{\mathrm{prox}_{\run}}

\newcommand{\epsi}{\epsilon}

\newcommand{\aux}{\zeta}

\newtheorem{lemma}{Lemma}

\newtheorem{remark}{Remark}

\newtheorem{assumption}{Assumption}

\title{Fast sparse optimization via adaptive shrinkage}

\author{Vito Cerone$~~~~$ Sophie M. Fosson\thanks{Department of Control and Computer Engineering, Politecnico di Torino, Italy (e-mail: sophie.fosson@polito.it). This paper is part of the project NODES which has received funding from the MUR – M4C2 1.5 of PNRR with grant agreement no. ECS00000036.}$~~~~$ Diego Regruto}
\date{}
\begin{document}
\maketitle

\begin{abstract}               
The need for fast sparse optimization is emerging, e.g., to deal with large-dimensional data-driven problems and to track time-varying systems. In the framework of linear sparse optimization, the iterative shrinkage-thresholding algorithm is a valuable method to solve Lasso, which is particularly appreciated for its ease of implementation. Nevertheless,   it converges slowly. In this paper, we develop a proximal method, based on logarithmic regularization, which turns out to be an iterative shrinkage-thresholding algorithm with adaptive shrinkage hyperparameter. This adaptivity substantially enhances the trajectory of the algorithm, in a way that yields faster convergence, while keeping the simplicity of the original method. Our contribution is twofold: on the one hand, we derive and analyze the proposed algorithm; on the other hand, we validate its fast convergence via numerical experiments and we discuss the performance with respect to state-of-the-art algorithms.
\end{abstract}



\section{Introduction}\label{sec:IN}
Sparse optimization consists in learning sparse models through the solution of suitable optimization problems. We call ``sparse'' those models that depend on a reduced  number of parameters, which is a desirable condition for several motivations, ranging from the decrease of  numerical complexity and memory footprint to the circumvention of overfitting; see, e.g., \cite{has15book,bru19book} for an overview. Recently, the attention on sparsity is increasing in system identification and in machine learning, on the one hand to select models that are physically interpretable, on the other hand, to train models that can be embedded in devices with limited resources, such as mobile applications; see, e.g., \cite{bru16,zha20,lou18} for different perspectives on this topic.

Nowadays, the optimization techniques to learn sparse models are quite mature, in particular in the context of linear systems. The Lasso problem, proposed in \cite{tib96}, is a popular, effective approach, that combines least squares and $\ell_1$ minimization to search sparse solutions via convex optimization. However, the development of fast algorithms for sparse optimization is still an open problem. In fact, in case of large-dimensional datasets, the convexity is not sufficient to guarantee solutions in a reasonable time. This aspect becomes critical, e.g., to identify time-varying or hybrid systems, where the parameters either evolve in time or switch among different unknown modes; in this framework, a prompt identification is necessary to track the dynamics; see, e.g., \cite{lau19,fox21}.

In the specific case of Lasso, several methods are proposed in the literature to minimize it efficiently. Among the iterative algorithms, the iterative shrinkage-thresholding algorithm (ISTA) proposed by \cite{dau04} is a simple proximal method, straightforward to implement even in a distributed context. Similarly, the alternating direction method of multipliers (ADMM) is low-complex and prone to distributed and parallelized computation; see \cite{boy10} for a complete overview. ISTA and ADMM share an R-linear convergence rate, as proven in \cite{bre08,hon17}, while in practice ADMM is usually faster. Methods to accelerate ISTA are proposed in the literature, such as FISTA, see \cite{bec09}, which improves the nonasymptotic global rate of convergence from $\mathcal{O}\left(t^{-1}\right)$ to $\mathcal{O}\left(t^{-2}\right)$.

In this paper, we analyze a different strategy for fast sparse optimization. We start from the development of a proximal gradient method for a Lasso-kind problem with logarithmic regularization, here denoted as Log-Lasso, and we obtain a variant of ISTA whose shrinkage hyperparameter is adaptive, i.e., it is updated at each iteration, based on the current estimate. This feature makes the algorithm significantly faster than ISTA. Throughout the paper, we name AD-ISTA the proposed adaptive ISTA.

More precisely, by starting from  the analysis in \cite{dah12}, we show that ISTA, FISTA and ADMM share a common behavior in the minimization of Lasso: in a first phase, they decrease the least squares term, while the $\ell_1$-norm substantially increases; in a second phase, they reduce the $\ell_1$-norm while keeping the least squares error almost constant. This trajectory yields a transient overestimation of the sparsity, which is time-consuming. In contrast, the proposed approach is not affected by this drawback.
In this work, we limit the investigation to linear systems and Lasso, while the proposed methodology might be extended to different optimization problems.

In summary, our contributions are as follows. Firstly, we illustrate how to apply a proximal method to Log-Lasso, which results into AD-ISTA, and we propose a straightforward analysis of its convergence; contextually, we show that AD-ISTA is equivalent to an ISTA with adaptive shrinkage parameter. Secondly, we illustrate via numerical experiments that AD-ISTA improves the trajectory with respect to classical algorithms, which results in a faster convergence. Moreover, by applying the approach of the fast ISTA (FISTA) proposed by \cite{bec09}, we propose an accelerated  version of AD-ISTA, called AD-FISTA. Finally, we also compare AD-ISTA and AD-FISTA to $\ell_1$ reweighting methods.

We organize the paper as follows. In Sec. \ref{sec:PS}, we state the problem. In Sec. \ref{sec:PA}, we derive and illustrate the proposed approach, while Sec. \ref{sec:why} discusses the expected faster behavior. This is validated via numerical experiments in Sec. \ref{sec:NR}, where we add to the comparison the related algorithms illustrated in Sec. \ref{sec:RA}; finally, we draw some conclusions.
\section{Problem Statement}\label{sec:PS}
Let us consider sparse optimization problems of the kind
\begin{equation}\label{p1}
\begin{split}
&\min_{x\in\R^n}\gun(x)+\runinitial(x)
\end{split}
\end{equation}
where $\gun :\R^n \mapsto \R^+ $ is a convex, smooth cost functional and
\begin{equation}\label{regu}
 \runinitial(x):= \sum_{i=1}^n \alpha_i r(x_i)
\end{equation}
is a sparsity promoting regularization with weights $\alpha=(\alpha_1,\dots,\alpha_n)$, $\alpha_i\geq 0$ for each $i=1,\dots,n$. Usually, $\runinitial(x)$ is not differentiable at $x=0$, and it may be non-convex.

To solve this class of composite non-smooth problems, where gradient descent is not feasible, we can exploit the proximal gradient method (PGM). PGM consists in iterating a Landweber step on $\gun$, i.e., a gradient descent step with constant stepsize $\tau>0$, and a  proximal mapping. More precisely, by defining $\lambda:=\tau\alpha\in\R^n$ and
\begin{equation}\label{prox_def}
\prox(z):=\argmin{x\in\R^n}\left(\run(x)+\frac{1}{2}\|x-z\|_2^2\right),
\end{equation}
PGM is as follows: for $t=0,1,2,\dots$
\begin{equation}\label{PGM}
\begin{split}
&x_{t+1}=\prox\left(x_t-\tau \nabla\gun(x_t)\right).
\end{split}
\end{equation}
The convergence of $\gun(x_t)+\runinitial(x_t)$, where $x_t$ is the sequence of iterates of PGM, is proven when $\run$ is convex, $\gun$ has Lipschitz continuous gradient, with constant $L$, and $\tau=\frac{2}{L}$. Specifically, if $\xstar\in\R^n$ is a minimizer, then $\gun(x_t)+\runinitial(x_t)$ converges to the  minimum $\gun(\xstar)+\runinitial(\xstar)$. In addition, if $\gun$ is strongly convex, the convergence of $x_t$ to the (unique) minimizer is also proven;  see, e.g., \cite{com11,cal14}.

Lasso is an instance of problem \eqref{p1}: given $A\in\R^{m,n}$ and a vector of measurements $y\in\R^m$, Lasso corresponds to $\gun(x)=\frac{1}{2}\|Ax-y\|_2^2$ and $\runinitial(x)=\alpha\|x\|_1$, with scalar $\alpha>0$. On the other hand, PGM applied to Lasso corresponds to ISTA. The convergence of ISTA can be proven in a peculiar way, developed by \cite{dau04}, that exploits the definition of a  surrogate  functional. With this approach, the sequence $x_t$ generated by ISTA is proven to converge  to a minimizer, even though $\gun(x)$ is not strongly convex. We notice that $\frac{1}{2}\|Ax-y\|_2^2$ is not strongly convex when $m<n$, which allows to efficiently use ISTA in compressed sensing; see, e.g., \cite{for10}.


The development of our approach starts by considering a Lasso with possible non-convex regularization $\runinitial$, i.e.,
\begin{equation}\label{p2}
\begin{split}
&\min_{x\in\R^n}\fun(x):=\frac{1}{2}\|Ax-y\|_2^2+\runinitial(x).
\end{split}
\end{equation}
In principle, we can apply PGM to solve \eqref{p2}. However, we have to cope with two critical points. On the one hand, the computation of $\prox$  may be not straightforward when $\run$ (or equivalently $\runinitial$) is non-convex. On the other hand, to prove the convergence of the algorithm is challenging.

In the literature, some attention is devoted to the application of PGM to problem \eqref{p2}.
%
%
In particular, \cite{bre15} provide conditions under which $\fun(x_t)$ converges in infinite-dimensional Hilbert spaces, by leveraging the results in \cite{att11}, and they analyze the specific case of $\ell_p$ regularization. Moreover, \cite{bay16} analyze the convergence of $x_t$ for \eqref{p2} when $\run(x)$ is weakly convex and  $A^TA$ is positive definite, which is not the case, e.g., of compressed sensing.

In the following, we focus on the case
\begin{equation}\label{r2}
r(x_i)=\log(|x_i|+\epsi)
\end{equation}
where $\epsi>0$ is a design hyperparameter that tunes the concavity of the function. We call Log-Lasso the problem \eqref{p2}-\eqref{regu}-\eqref{r2}. We consider Log-Lasso because log-regularization is known to be efficient, see, e.g., \cite{can08rew}; however, the proposed approach might be extended also to other classes of non-convex regularizers.

\section{Proposed approach}\label{sec:PA}
In this section, we develop the proposed approach by starting from the application of PGM to \eqref{p2}-\eqref{r2}. Firstly, we analyze the convergence of PGM applied to \eqref{p2}; then, we deal with the computation of the proximal operator of the log-regularizer.


Given $\fun(x)$ defined in \eqref{p2}, we introduce the surrogate functional
\begin{equation}\label{surrogate}
 \sun(x,\aux) = \fun(x)+\frac{1}{2\tau}\|x-\aux\|_2^2-\frac{1}{2}\|Ax-A\aux\|_2^2
\end{equation}
where $\aux\in\R^n$ is an auxiliary variable. We assume $\tau<\|A\|_2^{-2}$, so that
\begin{equation}\label{surrogate_pos}
 \frac{1}{\tau}\|x-\aux\|_2^2-\|Ax-A\aux\|_2^2\geq \left(\frac{1}{\tau}-\|A\|_2^2\right)\|x-\aux\|_2^2 \geq 0
\end{equation}
where the equality holds only for $x=\aux$. Then, we proceed by alternating minimization of \eqref{surrogate} with respect to $x$ and to $\aux$. According to \eqref{surrogate_pos},
\begin{equation}
 x = \argmin{\aux\in\R^n}\sun(x,\aux).
\end{equation}
On the other hand, we have
\begin{equation}\label{compute_proxy}
  \argmin{x\in\R^n}\sun(x,\aux) = \prox\left(\aux+\tau A^T(y-A\aux)\right).
\end{equation}
Then,
\begin{equation}\label{compute_proxy_dyn}
  x_{t+1}=\prox\left(x_t+\tau A^T(y-Ax_t)\right).
\end{equation}
We notice that
$\fun(x_{t})=\sun(x_{t},x_{t})\geq \sun(x_{t+1},x_{t})\geq \sun(x_{t+1},x_{t+1})= \fun(x_{t+1})\geq 0$, i.e.,  $\fun(x_{t})$ is a non-increasing, bounded sequence. This  yields the following result.
\begin{lemma}
Given the sequence $x_t$ generated by PGM applied to \eqref{p2}, $\fun(x_t)$ converges. Moreover, PGM applied to \eqref{p2} defines an asymptotically regular map, that is, $\lim\limits_{t\to\infty}\|x_{t+1}-x_t\|_2=0$.
\end{lemma}



In the following lemma, we specify how to compute $\prox(z)=\argmin{x\in\R^n}\run(x)+\frac{1}{2}\|x-z\|_2^2$ in case of log-regularizer \eqref{r2}, for any $z\in\R^n$. Since the problem is separable, we explicitly evaluate the minimum for each component $i=1,\dots,n$.
From now onwards, we tune the hyperparameters $\lambda=\tau\alpha$ and $\epsi$ such that
\begin{assumption}\label{ass_simple}
For each $i=1,\dots,n$, $\lambda_i<\epsi^2$.
\end{assumption}

\begin{lemma}\label{love} Under Assumption \ref{ass_simple}, given any $z_i\in\R$,
 \begin{equation}\label{boh}
\begin{split}
&\argmin{x_i\in\R}\mu(x_i):=\lambda_i \log(|x_i|+\epsi)+\frac{1}{2}(x_i-z_i)^2=\\
&=\left\{\begin{array}{ll}
                  z_i-\gamma_i(z_i)&~\text{ if } z_i> \frac{\lambda_i}{\epsi}\\
                  z_i+\gamma_i(z_i)&~\text{ if } z_i< -\frac{\lambda_i}{\epsi}\\
                  0&~\text{ if } z_i\in \left[-\frac{\lambda_i}{\epsi},\frac{\lambda_i}{\epsi}\right]
                 \end{array}
\right.
\end{split}
\end{equation}
where
\begin{equation}\label{def_gamma}
\gamma_i(z_i):=  \frac{|z_i|+\epsi -  \sqrt{(|z_i|+\epsi)^2-4\lambda_i}}{2}.
\end{equation}
\end{lemma}
{\textit{Proof. }}
Let us consider the case $x_i\in [0,+\infty)$.

For $x_i\in (0,+\infty)$, we have $\mu'(x_i)=\frac{\lambda_i}{x_i+\epsi}+x_i-z_i$, and $\mu''(x_i)=1-\frac{\lambda_i}{(x_i+\epsi)^2}>1-\frac{\lambda_i}{\epsi^2}$. Thus, $\mu''(x_i)>0$ under Assumption \ref{ass_simple}, i.e., $\mu(x_i)$ is strongly convex in $[0,+\infty)$. Therefore, if there exists a stationary point, it corresponds to the unique minimum; otherwise, the minimum is at $x_i=0$, since $\mu$ is strongly convex in $[0,+\infty)$ and $\lim\limits_{x_i\to+\infty}\mu(x_i)=+\infty$.

Let us compute the possible stationary points. We notice that $x_i+\epsi\neq 0$ because  $x_i\in[0,+\infty)$. Therefore,
\begin{equation}
\begin{split}
  \mu'(x_i)=0& ~~\Leftrightarrow \lambda_i+(x_i-z_i)(x_i+\epsi)=0\\
  &~~\Leftrightarrow x_i=\frac{z_i-\epsi\pm \sqrt{(z_i+\epsi)^2-4\lambda_i}}{2}.
\end{split}
\end{equation}
Then, we have to check whether these two solutions are consistent with the condition $x_i\geq 0$.

First of all, it is straightforward to evaluate that $(z_i+\epsi)^2-4\lambda_i\geq 0$ whenever $z_i\in\Omega :=(-\infty,-2\sqrt{\lambda_i}-\epsi)\cup(2\sqrt{\lambda_i}-\epsi,+\infty)$.
Afterwards,
\begin{equation*}
 \begin{split}
  &z_i-\epsi +\sqrt{(z_i+\epsi)^2-4\lambda_i}>0 \Leftrightarrow z_i\in\left(\frac{\lambda_i}{\epsi},+\infty\right)\subset \Omega; \\
  &z_i-\epsi -\sqrt{(z_i+\epsi)^2-4\lambda_i}<0 \text{ for any } z_i\in\Omega.
 \end{split}
\end{equation*}

We remark that $\left(\frac{\lambda_i}{\epsi},+\infty\right)\subset \Omega$ because $2\sqrt{\lambda_i}-\epsi<\frac{\lambda_i}{\epsi}$ under Assumption \ref{ass_simple}.

In conclusion, if $x_i\in[0,+\infty)$,  the minimizer is $ x_i=\frac{z_i-\epsi+ \sqrt{(z_i+\epsi)^2-4\lambda_i}}{2}$ if $z_i>\frac{\lambda_i}{\epsi}$; otherwise, if $z_i\leq \frac{\lambda_i}{\epsi}$, the minimizer is $x_i=0$.

We omit the computations for the case $x_i\in(-\infty, 0)$, which are based on the same ideas and yield symmetric conclusions.
To sum  up,
\begin{equation}\label{almost}
\begin{split}
&\argmin{x_i\in\R}\lambda_i \log(|x_i|+\epsi)+\frac{1}{2}(x_i-z_i)^2=\\
&=\left\{\begin{array}{ll}
                  \frac{z_i-\epsi +  \sqrt{(z_i+\epsi)^2-4\lambda_i}}{2}&~~\text{ if } z_i> \frac{\lambda_i}{\epsi}\\
                   \frac{z_i+\epsi- \sqrt{(z_i-\epsi)^2-4\lambda_i}}{2}&~~\text{ if } z_i< -\frac{\lambda_i}{\epsi}\\
                  0&~~\text{ otherwise.}
                 \end{array}
\right.
\end{split}
\end{equation}
From \eqref{almost}, we easily derive the thesis.
\QED
\begin{remark}
 The evaluation of the proximal operator is more complex for $\ell_p$. As one can see in \cite{zuo13,bre15}, no closed form is provided. This makes the choice of $\log$ regularization more affordable. On the other hand, the results in \cite[Lemma 3.3]{bre15} that define an exact formula for the proximal operation do not envisage the logarithmic case, as Assumption 3.2 in \cite{bre15} requires that $r'(x_i)\to \infty$ for $x_i \to 0$, which is not our case.
\end{remark}

Lemma \ref{love} provides an interesting interpretation of PGM applied to Log-Lasso. In fact, according to Lemma \ref{love}, we can formulate the PGM iteration as
\begin{equation}\label{adista}
\begin{split}
&z_t=x_t+\tau A^T(y-Ax_t),\\
&x_{t+1}=\soft_{\frac{\lambda}{\epsi}, \gamma(z_t)}\left(z_t\right)
\end{split}
\end{equation}
where $\gamma(z)=(\gamma_1(z_1),\dots,\gamma_n(z_n))$ is defined in \eqref{def_gamma} and $\soft_{\frac{\lambda}{\epsi}, \gamma(z)}(z)$ is a shrinkage-thresholding operator, defined componentwise as
\begin{equation}\label{softsoft}
 \soft_{\frac{\lambda_i}{\epsi}, \gamma_i(z_i)}(z_i):=\left\{\begin{array}{ll}
                  z_i-\gamma_i(z_i)&~~\text{ if } z_i> \frac{\lambda_i}{\epsi}\\
                  z_i+\gamma_i(z_i)&~~\text{ if } z_i< -\frac{\lambda_i}{\epsi}\\
                  0&~~\text{ otherwise.}
                 \end{array}
\right.
\end{equation}
In \eqref{softsoft}, $\gamma_i(z_i)$ is the shrinkage value, while $\frac{\lambda_i}{\epsi}$ is the threshold below which variables are set to zero.
We notice that, by using the notation of \eqref{softsoft}, we can write ISTA as
\begin{equation}\label{ista}
x_{t+1}=\soft_{\lambda,\lambda}\left(x_t+\tau A^T(y-Ax_t)\right)
\end{equation}
because in ISTA the shrinkage and thresholding hyperparameters are the same. In particular, in \eqref{ista}, $\lambda$ is time-invariant and constant for each component $i$.
%
In contrast, \eqref{adista} represents a generalization of ISTA, where shrinkage and thresholding hyperparameters, namely $\gamma(z)$ and $\frac{\lambda}{\epsi}$, are different. In particular, $\gamma(z)$ is time-varying, as it adapts to the current value of $x_t+\tau A^T(y-Ax_t)$ and it penalizes more the values closer to zero. This penalization is ``less democratic'' than the one of ISTA, and it causes less bias in the non-zero components of the solutions.
The proposed algorithm \eqref{adista} is an adaptive ISTA (AD-ISTA), because the shrinkage is adaptive.
\section{Why AD-ISTA is faster than ISTA?}\label{sec:why}
In this section, we discuss why AD-ISTA is expected to be faster than ISTA. The aim of this analysis is to illustrate the role played by the shrinkage hyperparameter in determining the trajectory of the algorithm. More rigorous proofs are left for extended work.

As to ISTA, the role of $\lambda$ in the speed of convergence is fundamental. In general, by increasing $\lambda$ we obtain a faster algorithm; on the other hand, a too large $\lambda$ would cause a substantial bias in the solution. Prior information on the solution can be used to set $\lambda$. As an extreme example, if the solution has support $\supp$, $|\supp|=k\ll n $, we assume $\lambda\in\R^n$ where $\lambda_i \notin \supp$ are very large, then the correct support is identified in one step, and the problem is immediately reduced to dimension $k\ll n $, which substantially reduces the number of iterations.

In \cite{dau08}, the Authors analyze the dynamics of ISTA for Lasso: ISTA first reduces the residual $\|Ax-y\|_2^2$ and contextually overestimates the $\ell_1$-norm; then, it corrects back the $\ell_1$-norm. This causes a ``long detour'' which yields slow convergence,  see \cite[Fig.1]{dau08}. For this motivation, the Authors propose to constrain the $\ell_1$-norm  within a given ball; however, this does not result in a clear acceleration of the method. 

The idea is further developed in \cite[Section 1.2]{dah09}, where the Authors propose to project the Landweber iteration $\ell_1$-balls with slowly increasing radius. In practice, they  implement the idea by proposing an ISTA with decreasing shrinkage-thresholding hyperparameter. This algorithm, called D-ISTA is further analyzed in \cite{dah12}. In particular,  it is proven to converge with R-linear rate under some conditions, e.g, by assuming that the hyperparameter decreases geometrically. However, choosing a suitable  decreasing hyperparameter is challenging.

In AD-ISTA we solve this issue, because the hyperparameter adapts to the magnitude of the gradient step over the previous estimate. This provides a larger shrinkage for small values in magnitude. In particular, if we start from the natural initial condition $x_0=0$, and if $\tau$ is small, in general the hyperparameter is larger in a first phase, which keeps the $\ell_1$-norm small, while as far as the components move away from zero, the shrinkage is smaller. Therefore, to some extent, we have a decreasing behavior, but only for those components that move far from zero. This adaptation is the key motivation that makes AD-ISTA more effective than ISTA.

We remark that in the literature the study of optimal hyperparameters for ISTA currently is an active topic. For example, a learned ISTA (LISTA) is developed in \cite{gre10}, and subsequently enhanced in, e.g., \cite{liu19,che21}. Since ISTA iterates a linear step and a non-linearity, the main structure is similar to a neural network, and techniques to learn the hyperparameters are studied in the mentioned papers. The main drawback of this approach is the time required for the training.

\section{Related algorithms}\label{sec:RA}
In this section, we present an accelerated version of AD-ISTA and we compare AD-ISTA to an $\ell_1$-reweighting method.

\subsection{AD-FISTA: a fast version of AD-ISTA}
Since AD-ISTA shares the same structure of ISTA, the fast version of ISTA proposed in \cite{bec09} and known as FISTA can be applied to it. Basically, FISTA exploits two previous iterates to compute the current estimate and it shares the improved convergence rate $\mathcal{O}\left(\frac{1}{t^2}\right)$, while keeping the low complexity of ISTA per iteration.

The application of FISTA approach to AD-ISTA, that we denote as AD-FISTA, is as follows: given $v_0=x_0\in\R^n,u_0 =1\in\R$, for any $t=0,1,2,\dots$,
\begin{equation}\label{adfista}
 \begin{split}
&z_t= v_t+\tau A^T(y-Av_t)\\
&x_{t+1}=\soft_{\frac{\lambda}{\epsi}, \gamma(z_t)}\left(z_t\right)\\
&u_{t+1}=\frac{1+\sqrt{1+4u_t^2}}{2}\\
&v_{t+1}=x_{t+1}+\frac{u_t-1}{u_{t+1}}(x_{t+1}-x_{t}).
 \end{split}
\end{equation}
For the same motivations illustrated in \cite{bec09}, we expect that AD-FISTA converges in less iterations than AD-ISTA.

\subsection{Comparison to $\ell_1$-reweighting ISTA}
In the context of sparse optimization,  $\ell_1$-reweighting techniques are popular to improve the accuracy of $\ell_1$ minimization methods. The key idea, proposed in \cite{can08rew} is to iterate the solution of a Basis Pursuit by updating the weight of the $\ell_1$-norm, with the final aim of penalizing less the larger components in magnitude, which is in line with the approach proposed in this work. The algorithm in \cite{can08rew} leverages the local minimization of a log-concave penalty through its linearization. This yields to weight the $\ell_1$-norm with the derivatives of the log-concave penalty. An ISTA-based variant of $\ell_1$ reweighting is proposed in \cite[Sec. III]{fox18}, and it is observed to be fast with respect to classic $\ell_1$-reweighting. In case of logarithmic penalty, this algorithm, here denoted as RW-ISTA, is as follows:
\begin{equation}\label{rwista}
\begin{split}
&(w_t)_i = \frac{1}{|(x_t)_i|+\epsi},~~~i=1,\dots,n\\
&x_{t+1}=\soft_{\lambda w_t,\lambda w_t}\left(x_t+\tau A^T(y-Ax_t)\right).
\end{split}
\end{equation}
Even though RW-ISTA originates from the $\ell_1$-reweighting framework, while AD-ISTA is obtained via proximal methods, the final structure of RW-ISTA has an adaptive  shrinkage-thresholding parameter as in AD-ISTA. This may explain the increased velocity of RW-ISTA observed in \cite{fox18}. On the other hand, in RW-ISTA, shrinkage and thresholding parameters are equal, while, as discussed above, in AD-ISTA they are different.
\section{Numerical results}\label{sec:NR}
In this section, we present some numerical results to validate the proposed method.

For our experiments, we consider a matrix $A\in\R^{m,n}$ with $m=500$ and $n=1000$, whose components are independently generated with Gaussian distribution $\mathcal{N}(0,\frac{1}{m})$. Given $y=A\xtrue+\eta$, where $\eta\in\R^m$ is an unknown random noise $\sim\mathcal{N}(0,10^{-2})$, we aim at estimating $\xtrue$, which has sparsity $k=10$, and non-zero components randomly generated with uniform distribution, with magnitude in $(1,2)$. To estimate $\xtrue$, we implement the proposed AD-ISTA and AD-FISTA, and we compare them to ISTA, FISTA, ADMM, and RW-ISTA. We set $\lambda=10^{-3}$ for Lasso, and $\lambda=4\times 10^{-4}$ and $\epsi = 10^{-2}$ for Log-Lasso. Finally, $\tau=\|A\|_2^{-2}$. The considered setting guarantees  that Lasso and Log-Lasso are successful, that is, by solving them we recover the correct support; see, e.g., \cite{fuc05} for theoretical guarantees. In particular, all the algorithms converge almost to the same solution. Instead, our goal is to analyze the convergence rate, in terms  of number of iterations. We specify that each iteration has comparable computational complexity for all the algorithms; therefore the number of iterations well represents the velocity of the algorithm.

\begin{figure}[h!]
\begin{center}
\includegraphics[width=0.7\columnwidth]{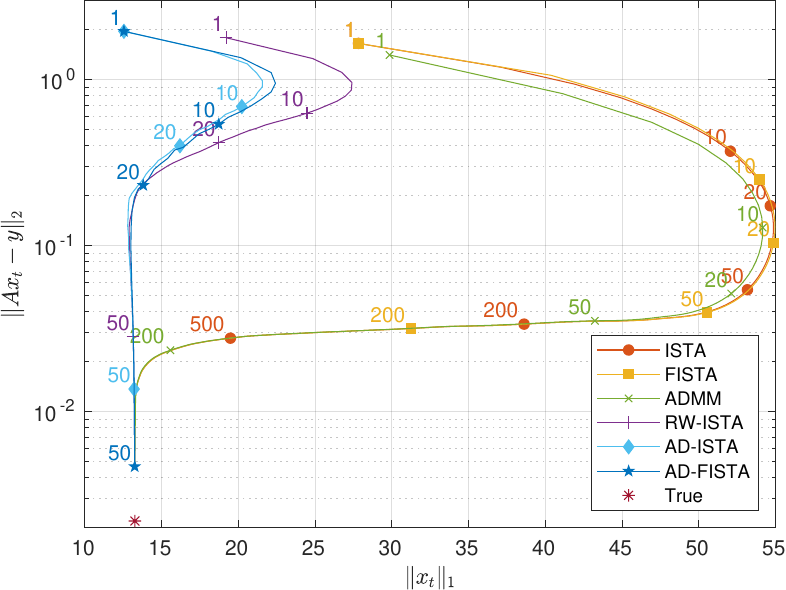}
\includegraphics[width=0.7\columnwidth]{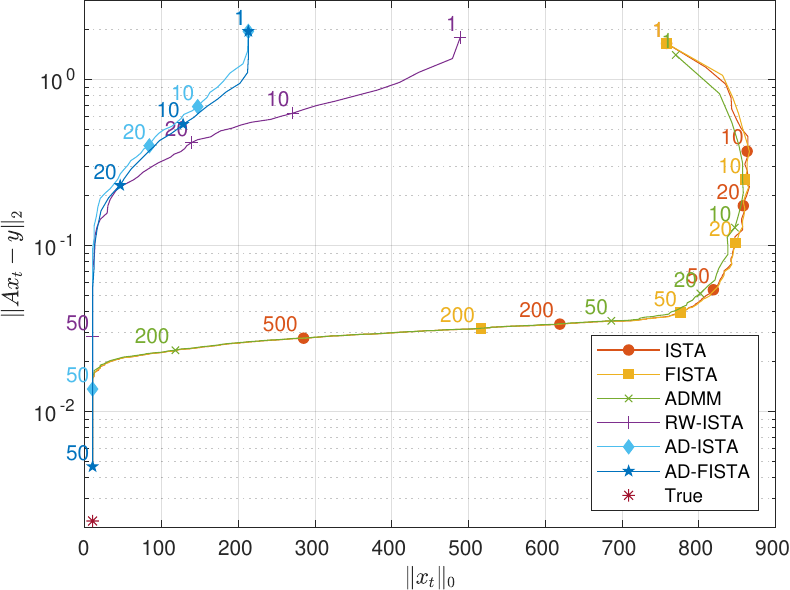}
\caption{Residual norm $\|Ax_t-y\|_2$ with respect to $\|x_t\|_1$ and $\|x_t\|_0$, respectively. The curve are parametrized with time. We label the iterations 1,10,20,50,200,500 to ease the comparison of the algorithms. ``True'' refers to the value of $\xtrue$, which is estimated by Lasso with a small bias on the non-zero components.}
\label{fig:1}
\end{center}
\end{figure}

In Fig. \ref{fig:1}, we can see the evolution of $\|Ax_t-y\|_2$ with respect to $\|x_t\|_1$ and $\|x_t\|_0$ in a single experiment. As discussed in previous sections, ISTA, FISTA, and ADMM are characterized by a transient overestimation of $\|x_t\|_1$. As a consequence, a similar behavior is observed for the sparsity $\|x_t\|_0$. In contrast, the proposed AD-ISTA and AD-FISTA keep $\|x_t\|_1$ smaller. As expected AD-FISTA is faster. Also RW-ISTA maintains a low  $\|x_t\|_1$, but less effectively than the proposed algorithms.

\begin{figure}[h!]
\begin{center}
\includegraphics[width=0.65\columnwidth]{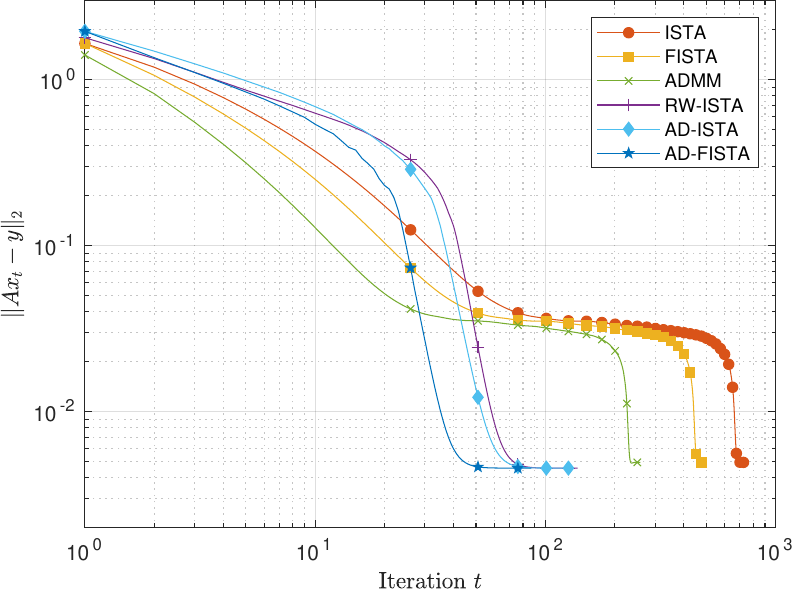}
\includegraphics[width=0.65\columnwidth]{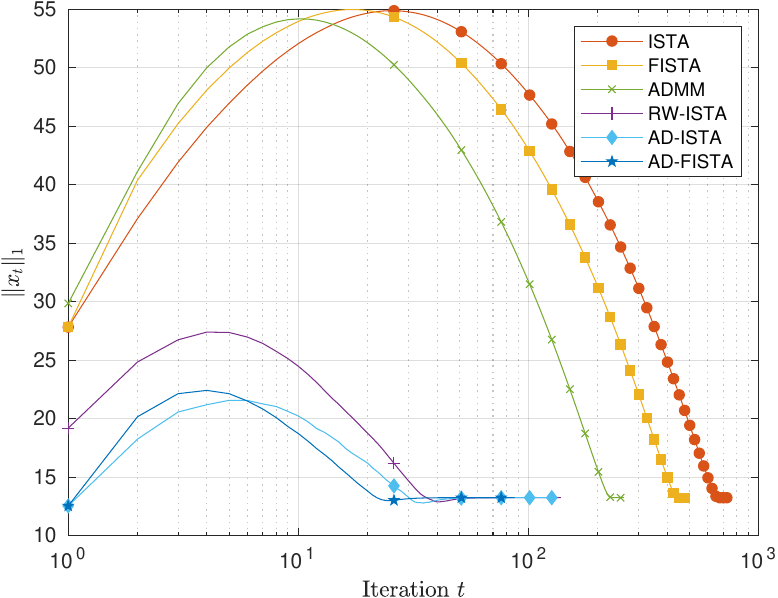}
\includegraphics[width=0.65\columnwidth]{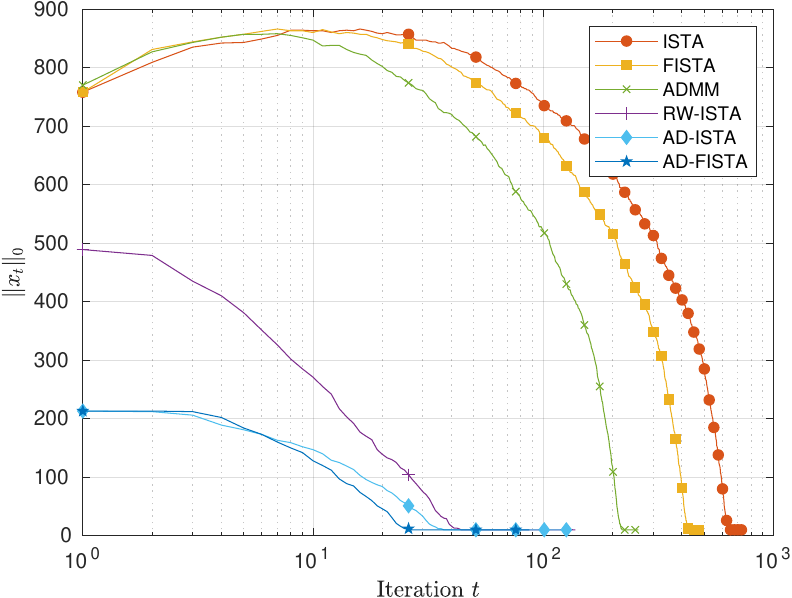}
\caption{The time evolution of $\|Ax_t-y\|_2$, $\|x_t\|_1$ and $\|x_t\|_0$.}
\label{fig:2}
\end{center}
\end{figure}

In Fig. \ref{fig:2}, we depict the time evolution of $\|Ax_t-y\|_2$, $\|x_t\|_1$ and $\|x_t\|_0$, for the same experiment. We can see that the number of iterations required by AD-ISTA and AD-FISTA is substantially lower than the one of ISTA, FISTA and ADMM.

\begin{table}[h]
\renewcommand{\arraystretch}{1.4}
\centering
    \begin{tabular}{c | c c c }
     & Number & of & iterations\\
    \hline
    Algorithm & Mean & Min & Max \\
    \hline

  ISTA &895.44& 703& 1085\\
  FISTA &595.94& 467&722\\
  ADMM & 318.49&255&378\\
  RW-ISTA&147.47&126&173\\
  AD-ISTA&138.34&119&162\\
  AD-FISTA&90.64&78&107\\
  \end{tabular}
    \caption{Number of iterations to converge over 100 random runs.}
    \label{tab:1}
\end{table}
In Table \ref{tab:1}, we collect some statistics on the number of iterations over 100 random runs. We notice that the maximum of iterations required by AD-ISTA and AD-FISTA is smaller than the minimum of iterations required by ISTA, FISTA and ADMM.

\section{Conclusions}\label{sec:C}
In this work, we propose and analyze AD-ISTA, a variant of ISTA developed by applying the proximal gradient method to Log-Lasso. AD-ISTA converges in less iterations with respect to ISTA, FISTA and ADMM, thanks to an adaptive shrinkage hyperparameter, that limits the increase of the $\ell_1$-norm during the first phase. Moreover, by applying  the principles of FISTA, we also propose the accelerated version AD-FISTA. Through numerical experiments, we verify that AD-ISTA is faster than the state-of-the-art algorithms for Lasso and that we obtain a further acceleration with AD-FISTA. Possible extensions of this work include the rigorous proof of the convergence rate and the generalization to sparse optimization problems different from Lasso.
%
\bibliographystyle{IEEEtran}
\bibliography{refs}             
\end{document}